\documentclass[12pt]{article}
\usepackage{amssymb}
\usepackage{amsmath}
\usepackage{latexsym}
\usepackage{mathrsfs}

\numberwithin{equation}{section}
\newtheorem{thm}[equation]{Theorem}

\newtheorem{rem}[equation]{Remark}

\newtheorem{cor}[equation]{Corollary}

\providecommand{\abs}[1]{\left\lvert#1\right\rvert}

\title{Viewing the Steklov eigenvalues of the Laplace operator as critical Neumann eigenvalues}

\author{Pier Domenico Lamberti\footnote{Corresponding author. The results discussed in this paper have been announced by the second author at the  9th
ISAAC Congress, Krak\'{o}w 2013.}\,\,   and Luigi Provenzano}

\date{\ }

\begin{document}

\newcommand{\rea}{\mathbb{R}}

\maketitle

%
%
%

\noindent
{\bf Abstract:}
We consider the Steklov eigenvalues of the Laplace operator as limiting Neumann eigenvalues in a problem of boundary mass concentration.  We discuss the asymptotic behavior of the Neumann eigenvalues in a ball and we deduce 
that the Steklov eigenvalues minimize the  Neumann eigenvalues. Moreover, we study the dependence of the eigenvalues of the Steklov problem upon 
perturbation of the mass density  and show that the Steklov eigenvalues violates a maximum principle in spectral optimization problems.

\vspace{11pt}


\vspace{6pt}

\section{Introduction}

Let $\Omega$ be a bounded domain (i.e. a bounded connected open set)  of class $C^2$ in $\mathbb R^N$, $N\geq 2$. We consider the  Steklov 
eigenvalue problem for the Laplace operator 
\begin{equation}\label{Ste}
\left\{\begin{array}{ll}
\Delta u =0 ,\ \ & {\rm in}\  \Omega,\\
\frac{\partial u}{\partial \nu }=\lambda\rho u ,\ \ & {\rm on}\  \partial \Omega,
\end{array}\right.
\end{equation}
in the unknowns $\lambda$ (the eigenvalue) and $u$ (the eigenfunction). Here $\rho $ denotes a positive function on $\partial \Omega$ bounded away from zero and infinity and $\nu$  the unit outer normal to $\partial\Omega$. 

Keeping in mind important problems in linear elasticity (see e.g. Cou\-rant and Hilbert \cite{cohi}), we shall think of the weight $\rho$ as a mass density.
In fact, for $N=2$ problem (\ref{Ste}) arises for example in the study of the vibration modes of a free elastic membrane the total mass of which is concentrated at the boundary. Note that the total mass is given by $\int _{\partial \Omega}\rho d\sigma $. This mass concentration phenomenon can be 
described as follows. 

For any $\epsilon >0$ sufficiently small, we consider the $\epsilon$-neighborhood of the boundary 
$
\Omega_{\epsilon}=\{x\in\Omega:d(x,\partial\Omega)<\epsilon\}
$
and for a fixed $M>0$ we define a function $\rho_{\epsilon}$ in the whole of $\Omega$ as follows

\begin{equation}\label{dens}
\rho_{\epsilon}(x)=\left\{
\begin{array}{ll}
\epsilon,& {\rm if\ }x\in \Omega\setminus\overline\Omega_{\epsilon},\\
\frac{M-\epsilon\abs{\Omega\setminus\overline\Omega_\epsilon}}{\abs{\Omega_{\epsilon}}},& {\rm if\ }x\in \Omega_\epsilon.
\end{array}
\right.
\end{equation}
Note that for any $x\in\Omega $  we have $\rho_{\epsilon}(x)\to 0$ as $\epsilon \to 0$, and   $\int_{\Omega}\rho_{\epsilon}dx=M$ for all $\epsilon>0$. Then we consider the following eigenvalue problem 
for the Laplace operator with Neumann boundary conditions
\begin{equation}\label{Neu}
\left\{\begin{array}{ll}
-\Delta u =\lambda \rho_{\epsilon} u,\ \ & {\rm in}\  \Omega,\\
\frac{\partial u}{\partial \nu }=0 ,\ \ & {\rm on}\  \partial \Omega. 
\end{array}\right.
\end{equation}
We recall  that for $N=2$ problem (\ref{Neu}) provides the vibration modes of a free elastic membrane with mass density $\rho_{\epsilon }$ and total mass
$M$. 
It is not difficult to prove that  the eigenvalues and eigenfunctions of problem (\ref{Neu}) converge as $\epsilon$ goes to zero to the eigenvalues and eigenfunctions of problem (\ref{Ste}) with  $\rho =\frac{M}{\abs{\partial\Omega}}$. Thus the Steklov problem can be considered as a limiting Neumann problem.  We refer to \cite[Arrieta, Jim\'{e}nez-Casas, Rodr\'{\i}guez-Bernal]{arr} for a general approach to this type of problems. 

The aim of this paper is to highlight a few properties of the Steklov problem which, compared to the  Neumann problem, reveals a critical nature.

First,  we study the asymptotic behavior of the eigenvalues of problem (\ref{Neu}) as $\epsilon \to 0$, when $\Omega$ is a ball. We prove that such eigenvalues are differentiable with  respect to $\epsilon \geq 0$ and establish formulas for the first order derivatives at $\epsilon =0$, see Theorem~\ref{deriv}. It turns our that such derivatives are positive, hence the Steklov eigenvalues minimize the Neumann eigenvalues of problem (\ref{Neu}) for $\epsilon$ sufficiently small, see Remark~\ref{monrem}.  

Second, we consider the problem of optimal mass distributions for problem (\ref{Ste}) under the condition that  that the total mass is fixed. This  problem has been largely investigated in the case of Dirichlet boundary conditions, see e.g. Henrot~\cite{he} for references. As for Steklov boundary conditions, we quote
the classical paper by Bandle and Hersch~\cite{banher}.  

By following the approach developed in \cite{lam}, we prove that simple eigenvalues and the symmetric functions of the multiple eigenvalues of (\ref{Ste}) depend real analytically on $\rho$ and we characterize the corresponding critical mass densities under mass constraint. See Theorem~\ref{sym} and Corollary~\ref{crit}. Again, the Steklov problem exhibits a critical behavior and violates the maxi\-mum principle discussed in \cite{lampro}  for general elliptic operators of arbitrary order subject to homogeneous boundary conditions of Dirichlet, Neumann and intermediate type for  which critical mass densities do not exist. Indeed, it turns out that if $\Omega$ is a ball then the constant function  is a critical mass density for the Steklov problem (\ref{Ste}), see  Corollary~\ref{critste}, Remark~\ref{rembandle} and Theorem~\ref{teobandle}.

\section{Asymptotic behavior of Neumann eigenvalues}

Given a bounded domain $\Omega$ in ${\mathbb{R}}^N$  of class $C^2$  and $M>0$  we denote by 
$\lambda_j$, $j\in\mathbb N$,  the eigenvalues of problem (\ref{Ste}) corresponding to the constant surface density $\rho=\frac{M}{\abs{\partial\Omega}}$. Similarly,   for $\epsilon >0$ sufficiently small, we denote by 
$\lambda_j(\epsilon )$, $j\in\mathbb N$, the eigenvalues of problem (\ref{Neu}).  Note that in this paper we always assume that  $N\geq 2$. Moreover, by ${\mathbb{N}}$ we denote the set of natural numbers including zero, hence   $\lambda_0(\epsilon)=\lambda_0=0$ for all $\epsilon >0$.

As is well-known,  by the Min-Max Principle we get the following  
variational characterization of the two sequences of eigenvalues:

\begin{eqnarray*}
&&\lambda_{j}(\epsilon )=\inf_{\substack{E\subset  H^{1}(\Omega)  \\{\rm dim}E=j+1}}\sup_{0\ne u\in E}\frac{\int_\Omega\abs{\nabla  u}^2 dx}{\int_{\Omega} u^{2}\rho_{\epsilon}\,dx},\ \ \ \mathcal8 j\in\mathbb N,\\
&&\lambda_{j}=\inf_{\substack{ E\subset  H^{1}(\Omega)\\{\rm dim} E=j+1}}\sup_{\substack{ u\in  E \\ {\rm Tr}\, u\ne 0 }}\frac{\int_\Omega\abs{\nabla  u}^2 dx}{\int_{\partial\Omega}  ({\rm Tr}\, u) ^{2}\frac{M}{|\partial \Omega |}\,d\sigma},\ \ \ \mathcal8 j\in\mathbb N.
\end{eqnarray*}
Here $H^1(\Omega )$ denotes the standard Sobolev space of real-valued functions in $L^2(\Omega)$ with weak derivatives up to first order in $L^2(\Omega)$ and ${\rm Tr}\, u$ denotes the trace in $\partial \Omega$ of a function $u\in H^1(\Omega)$ . 
We note that, for each fixed $u\in H^{1}(\Omega)$ we have
\begin{equation}\label{limray}
\lim_{\epsilon\rightarrow 0}\frac{\int_\Omega\abs{\nabla  u}^2 dx}{\int_{\Omega}  u^{2}\rho_{\epsilon}\,dx}=\frac{\int_\Omega\abs{\nabla  u}^2 dx}{\int_{\partial\Omega}({\rm Tr}\, u)^{2}\frac{M}{|\partial \Omega|}  \,d\sigma}.
\end{equation}
By looking at (\ref{limray}) one could expect the spectral convergence of  the Neumann problems  under consideration  to the Steklov problem. In fact the following statement holds.

\begin{thm}\label{arrietaanibal} If $\Omega $ is a bounded domain in ${\mathbb{R}}^N$ of class $C^2$ then  $\lim_{\epsilon\rightarrow 0}\lambda_j(\epsilon )=\lambda_j$ for all $j\in\mathbb N$.
\end{thm}

\noindent This theorem can be proved directly  by using the notion of compact convergence for the resolvent operators  but can also be obtained as a consequence of the more general results proved  in \cite[Arrieta, Jim\'{e}nez-Casas, Rodr\'{\i}guez-Bernal]{arr}. 

By Theorem \ref{arrietaanibal}, it follows that the function $\lambda_j(\cdot )$ can be extended with continuity at $\epsilon =0$ by setting   $\lambda_j(0)=\lambda_j$ for all $j\in {\mathbb{N}}$. This will be understood in the sequel.
If $\Omega$ is a ball then we are able to establish the asymptotic behavior of $\lambda_j(\epsilon )$ as $\epsilon \to 0$. Indeed, we can prove that
$\lambda_j(\epsilon )$ is differentiable with respect to $\epsilon$ and compute the derivative $\lambda_j'(0)$ at $\epsilon =0$. 

\begin{thm}\label{deriv}
If $\Omega$ is  the unit ball in $\mathbb R^N$ then $\lambda_j(\epsilon )$ is differentiable for any $\epsilon \geq 0$ sufficiently small and 
$$
\lambda_j'(0)=\frac{2M\lambda_j^2(0)}{3N|\Omega |}+\frac{2\lambda^2_j(0)|\Omega|}{2M\lambda_j(0)+N^2|\Omega|}.
$$
\end{thm}

\noindent The proof of this theorem relies on the use of Bessel functions which allow to recast the Neumann eigenvalue problem in the form of an equation 
$F(\lambda , \epsilon)=0$ in the unknowns $\lambda ,\epsilon$.  Then, after some preparatory work, it is possible to apply   the Implicit Function Theorem  and conclude. We note that, despite the idea of the proof is rather simple and used also in other contexts (see e.g. \cite{lape}), this method requires standard but lengthy computations, suitable Taylor's expansions and estimates on the corresponding remainders, as well as recursive formulas for the
cross-products of Bessel functions and their derivatives. We refer to \cite{proz} for details.

\begin{rem}\label{monrem} By Theorem \ref{deriv} it follows that for $\epsilon > 0$ sufficiently small the functions $\epsilon \mapsto \lambda_j(\epsilon )$ are strictly  increasing. 
In particular, it follows that for all $\epsilon >0$ sufficiently small, we have that  
$
\lambda_j(0)< \lambda_j(\epsilon)$.

It is interesting to compare our result with the monotonicity result by Ni and Wang~\cite{niwa} who have proved that if $\Omega$  is the unit disk in the plane  then the first positive eigenvalue
of the Neumann Laplacian in $\Omega_{\epsilon}$, i.e. the first positive eigenvalue of the problem 
\begin{equation}\label{Neuniwa}
\left\{\begin{array}{ll}
-\Delta u =\lambda u,\ \ & {\rm in}\  \Omega_{\epsilon },\\
\frac{\partial u}{\partial \nu }=0 ,\ \ & {\rm on}\  \partial \Omega_{\epsilon},
\end{array}\right.
\end{equation}
is a strictly increasing function of $\epsilon >0$. 
\end{rem}

\section{Existence of critical mass densities for the Steklov problem} \label{subcrit}

Given a bounded domain $\Omega$ in ${\mathbb{R}}^N$ of class $C^2$,  we denote by ${\mathcal R}$ the subset of $L^{\infty}(\partial\Omega)$ of those functions $\rho\in L^{\infty}(\partial\Omega)$ such that ${\rm ess }\inf _{\partial\Omega }\rho >0$. For any $\rho \in {\mathcal{R}}$, we denote by $\lambda_j[\rho ]$, $j\in {\mathbb{N}}$, the eigenvalues of problem 
(\ref{Ste}). By classical results in perturbation theory, one can prove that $\lambda_j[\rho]$ depends real-analytically on $\rho$ as long as $\rho $ is such that
$\lambda_j[\rho ]$ is a simple eigenvalue. This is no longer true if the multiplicity of $\lambda_j[\rho]$  varies. As it was pointed out in \cite{lam, lala2004}, in the case of multiple
eigenvalues, analyticity can be proved for the symmetric functions of the eigenvalues.  Namely, given a finite set of indexes  $F\subset {\mathbb{N}}$, one can consider 
the symmetric functions of the eigenvalues with indexes in $F$
\begin{eqnarray}
\label{sym1}
\Lambda_{F,h}[\rho ]=\sum_{ \substack{ j_1,\dots ,j_h\in F\\ j_1<\dots <j_h} }
\lambda_{j_1}[\rho ]\cdots \lambda_{j_h}[\rho ],\ \ \ h=1,\dots , |F|  \nonumber
\end{eqnarray}
and prove that such functions are real-analytic on    
\begin{equation}
{\mathcal { R}}[F]\equiv \left\{\rho\in {\mathcal { R}}:\
\lambda_j[\rho ]\ne \lambda_l[\rho ],\ \forall\  j\in F,\,   l\in \mathbb{N}\setminus F
\right\}.
\end{equation}
In fact, we can prove the following theorem where in order to establish formulas for the Frech\'{e}t differentials, we find it convenient to set
\begin{eqnarray*}
\Theta [F]\equiv \left\{\rho\in {\mathcal { R}}[F]:\ \lambda_{j_1}[\rho ]
=\lambda_{j_2}[\rho ],\, \
\forall\ j_1,j_2\in F  \right\} .
\end{eqnarray*}

\begin{thm}
\label{sym}
Let $\Omega$ be a bounded domain in $\mathbb R^N$ of class $C^2$ and $F$ a finite subset of ${\mathbb{N}}$. Then  ${{\mathcal { R}}}[F]$ is an open set in $L^{\infty}(\partial\Omega )$ and the functions
$\Lambda_{F,h}$ 
are real-analytic in ${{\mathcal { R}}}[F]$. Moreover, if $F=\cup_{k=1}^nF_k$ and
$\rho\in \cap_{k=1}^n\Theta [F_k]$ is such that for each $k=1,\dots , n$ the eigenvalues $\lambda_j[\rho ]$ assume the common value $\lambda_{F_k}[\rho ]$
for all $j\in F_k$, then the differentials of the functions $\Lambda_{F,h}$
at the point $\rho$ are given by the formula

\begin{eqnarray}
\label{sym2}
d\Lambda_{F,h}[\rho][\dot{\rho}] =-\sum_{k=1}^n
 c_k
\sum_{l\in F_k}
\int_{\partial\Omega}({\rm Tr}\, u_l)^2\dot{\rho}d\sigma\, ,
\end{eqnarray}
for all $\dot{\rho}\in L^{\infty }(\partial\Omega)$, where
$$ c_k=
\sum_{\substack{0\le h_1\le |F_1|\\ \dots\dots \\ 0\le h_n\le  |F_n|\\ h_1+\dots +h_n=h }}
{ |F_k|-1 \choose h_k-1 }\lambda_{F_k}^{h_k}[\rho]  \prod_{\substack{j=1\\ j\ne k}}^n { |F_j| \choose h_j }\lambda_{F_j}^{h_j}[\rho] ,
$$
and for each $k=1,\dots , n$, $\{  u_l\}_{l\in F_k}$ is a  basis of  the eigenspace of  $\lambda_{F_k}[\rho ]$ normalized by the condition $\int_{\partial \Omega}{\rm Tr}\, u_i{\rm Tr}\, u_j\rho d\sigma=\delta_{ij}$ for all $i,j\in F_{k}$.
\end{thm}

The proof of this theorem follows the lines of the corresponding result proved in \cite{lampro} for general elliptic operators subject to homogeneous boundary conditions of Dirichlet, Neumann and intermediate type. In the same spirit of \cite{lampro}, we can use formula (\ref{sym2}) in order to investigate the existence of critical mass densities  for the eigenvalues of the Steklov problem subject to mass constraint. 
We note that a typical optimization problem in the analysis of composite materials consists in finding mass densities $\rho$, with given total mass, which minimize a cost functional $F[\rho ]$ associated with the solutions of suitable partial differential equations depending on $\rho$.   Namely, in the case of Steklov boundary conditions one can consider the following problems   
$$
\min_{\int_{\partial \Omega}\rho d\sigma ={\rm const.}}F[\rho ]\ \ \ {\rm or}\ \ \ \max_{\int_{\partial \Omega }\rho d\sigma ={\rm const.}}F[\rho ].
$$
More in general, setting $M[\rho ]=\int_{\partial \Omega}\rho d\sigma $ one can consider the problem of finding critical mass densities $\rho$ under mass constraint, i.e. mass densities $\rho$ which satisfy the condition
$
{\rm Ker }dM[\rho ]\subset {\rm Ker }dF[\rho ]. 
$
As in \cite{lampro} we can give a characterization of critical mass densities which immediately follows by formula (\ref{sym2}) combined with the Lagrange Multipliers Theorem.

\begin{cor}\label{crit}
Let all  assumptions of Theorem \ref{sym} hold. Then, $\rho\in {\mathcal {R}}$ is a critical mass density for $\Lambda_{F,h}$  for some  $h=1,...,\abs{F}$, subject to  mass constraint if and only if there exists $c\geq 0$  such that
\begin{eqnarray}\label{overd}
 \sum_{k=1}^nc_k\sum_{l\in F_k}({\rm Tr}\, u_l)^2  =c, \ \ {\rm a.e.\ on }\ \partial \Omega.
\end{eqnarray}
\end{cor} 

The analysis carried out in \cite{lampro} has pointed out that for a large class of non-negative elliptic operators subject to homogeneous boundary conditions of intermediate type (including the  case of Dirichlet boundary conditions), there are no critical mass densities for simple eigenvalues and the symmetric functions of multiple eigenvalues. For example,  in the case of Dirichlet or Neumann boundary conditions,  
(\ref{overd}) has to be replaced by 
\begin{eqnarray}\label{overdbis}
 \sum_{k=1}^nc_k\sum_{l\in F_k}u_l^2  =c, \ \ {\rm a.e.\ in }\  \Omega.
\end{eqnarray}
which is clearly not satisfied in the Dirichlet case. As for Neumann boundary conditions the same non existence result can be easily proved for simple eigenvalues in which case  only a summand appears in (\ref{overdbis}). The situation is not completely clear for multiple eigenvalues. Under suitable regularity assumptions on the eigenfunctions $u_1$ and $u_2$ associated with  the same Neumann eigenvalue $\lambda$ one can prove that the condition $u_1^2+u_2^2=c$ in $\Omega$ implies that $\lambda=0$, but the proof in the case of multiplicities higher than two seems not straightforward. However, well-known explicit formulas for the eigenfunctions of the Neumann Laplacian in the ball clearly show that condition (\ref{overdbis}) is not satisfied, hence {\it no critical mass densities exist for the Neumann Laplacian in the ball}.   
In the case of Steklov boundary conditions the situation is much different. Indeed, if $\Omega$ is a ball then a critical mass density exists. 

\begin{cor}\label{critste}
Let $\Omega$ be the unit ball in $\mathbb R^N$,  $M>0$ and $F\subset {\mathbb{N}}$ be a finite set such that  the constant mass density $\rho=M/|\partial \Omega |$ belongs to ${\mathcal{R}}[F]$. Then  $\rho=M/|\partial \Omega |$ is  critical  for $\Lambda_{F,h}$ for all $h=1,...,|F|$ under the constraint $M[\rho ]=M$.
\end{cor}

\noindent The proof can be carried out as in \cite{lalacri}. Namely,  assume that  $\lambda $ is an eigenvalue of problem (\ref{Ste}) with multiplicity $m$ and consider a
basis  $u_1, \dots , u_m$ of the corresponding eigenspace. Assume that this basis is orthonormal in $L^2(\partial \Omega)$ with respect to the scalar product defined by $\int_{\partial \Omega}{\rm Tr}\, u{\rm Tr}\, v\rho d\sigma$. Then for any isometry $R$ in ${\mathbb{R}}^N$ also $u_1\circ R, \dots , u_m\circ R$ is an  orthonormal 
basis of the same eigenspace, hence $\sum_{i=1}^mu_i^2= \sum_{i=1}^mu_i^2\circ R$. It follows  that  $\sum_{i=1}^mu_i^2$ is constant on $\partial \Omega$.

\begin{rem}  \label{rembandle}
It is interesting to compare Corollary~\ref{critste} with a classical result proved by Bandle and Hersch~\cite{bandle} in the case of a class of symmetric planar domains. For the convenience of the reader  we formulate such result assuming directly that  $\Omega $ is the  unit disk in ${\mathbb{R}}^2$ centered at zero.
For any $n\in {\mathbb{N}}$ we set
$$
{\mathcal{R}}_{n}=\{\rho \in {\mathcal{R}}:\ \rho (e^{2\pi i/n}z)=\rho (z),\  \forall \ z\in \partial \Omega \},
$$
where the use of the complex variable $z$ is clearly understood. Then we have the following result
\begin{thm}[Bandle and Hersch]\label{teobandle}
Let $\Omega $ be the unit disk in ${\mathbb{R}}^2$ centered at zero, $M>0$, $n\in {\mathbb{N}}$. Then
$$
\lambda_j[\rho]\le \lambda_j\left[\frac{M}{2\pi}\right]
$$
for all all $j=0,\dots , n$ and $\rho \in {\mathcal{R}}_{ n}$ such that $M[\rho]=M$ . Equality holds only if $\rho =M/2\pi$.
\end{thm}
\noindent Thus in the case of a ball in ${\mathbb{R}}^2$ the constant mass density is in fact a maximizer among all mass densities satisfying the symmetry condition above. We refer to Bandle~\cite{bandle} for further discussions.  
\end{rem}

{\bf Acknowledgments}: We acknowledge financial support  by the research pro\-ject  ``Singular perturbation problems for differential operators", Progetto di Ateneo
of the University of Padova. 
\\

\noindent {\small
Pier Domenico Lamberti and Luigi Provenzano\\
Dipartimento di Matematica\\
Universit\`{a} degli Studi di Padova\\
Via Trieste,  63\\
35126 Padova\\
Italy\\
e-mail:	lamberti@math.unipd.it \\
e-mail: proz@math.unipd.it

}
\end{document}